\documentclass[11pt,bezier]{article}
\usepackage{amsmath,amssymb,amsfonts,euscript}
\usepackage[usenames]{color}
\newpage

\newcommand{\be}{\begin{equation}}
\newcommand{\ee}{\end{equation}}
\newtheorem{theorem}{Theorem}[section]
\newtheorem{proof}{Proof}[section]

\newtheorem{corollary}{Corollary}[section]

\title{\bf\Large Parameter Estimation of Type-II Hybrid Censored Weighted Exponential Distribution}
\vspace{-1cm}
\author
{A.Kohansal., S.Rezakhah \footnote{Faculty of Mathematics and Computer Science, Amirkabir University of Technology,
Tehran, Iran. Email: rezakhah@aut.ac.ir, ak\_kohansal@aut.ac.ir}}
\begin{document}
\maketitle

\begin{abstract}
A hybrid censoring scheme is a mixture of Type-I and Type-II censoring schemes. We study the estimation of parameters of weighted exponential distribution based on Type-II hybrid censored data. By applying EM algorithm, maximum likelihood estimators are evaluated. Also using Fisher infirmation matrix asymptotic confidence intervals are provided. By applying Markov Chain Monte Carlo techniques Bayes estimators, and corresponding highest posterior density confidence intervals of parameters are obtained. Monte Carlo simulations to compare the performances of the different methods is performed and one data set is analyzed for illustrative purposes.

\end{abstract}

\section{Introduction}
\par
Type-I and Type-II censoring schemes are two most popular censoring schemes which are used in practice. They can be briefly described as follows. Suppose $n$ units are put on a life test. In Type-I censoring, the test is terminated when a pre-determined time, $T$, on test has been reached, and failures after time
$T$ are not observed. In Type-II censoring, the test is terminated when a pre-chosen number, $R$, out of $n$ items has failed. It is also assumed that the failed items are not replaced. So, in Type-I censoring scheme, the number of failures is random and in Type-II censoring scheme, the experimental time is random.\par
A hybrid censoring scheme is a mixture of Type-I and Type-II censoring schemes and it can be described as follows. Suppose $n$ identical units are put to test. The test is finished when a pre-selected number $R$ out of $n$ items are failed, or when a pre-determined time $T$ on the test has been obtained. From
now on, we call this Type-I hybrid censoring scheme and this scheme has been used as a reliability acceptance test
in \cite{mil}. This censoring scheme was introduced by Epstin \cite{eps}, he also studied the life testing data under the assumption of exponential distribution with mean life $\theta$. Epstein \cite{eps} proposed two-sided confidence intervals for $\theta$ without any formal proof. Fairbanks et al. \cite{fai} moderated  partly the proposition of Epstein \cite{eps} and suggested a simple set of confidence intervals. Chen and Bhattacharya \cite{che} earned the exact distribution of the conditional maximum likelihood estimator (MLE) of $\theta$ and implied a one-sided confidence interval. Childs et al. \cite{chil} proposed some simplifications of the exact distribution. From the Bayesian point of view, Drapper and Guttmann \cite{dra} studied the same problem, and reached a two-sided credible interval of the mean lifetime based on the gamma prior. Comparison of the different methods using Monte Carlo simulations, can be found in Gupta and Kundu \cite{gup3}. For some related work, one may refer to Ebrahimi \cite{ebra1,ebra2}, Jeong et al. \cite{jeo}, Childs et al. \cite{chil}, Kundu \cite{kun1}, Banerjee and Kundu \cite{ban}, Kundu and Pradhan \cite{kun2}, Dube et al. \cite{dub} and the references cited there.
\par
One of the disadvantages of Type-I hybrid censoring scheme is that there may be very few failures occurring up to the pre-fixed time $T$. Because of this, Childs et al. \cite{chil} proposed a new hybrid censoring scheme known as Type-II hybrid censoring scheme which can be described as follows. Put $n$ identical items on test, and then stop the experiment at the random time $T^*=\mbox{max}\{x_{R:n},T\}$ , where $R$, and $T$ are prefixed numbers and $x_{R:n}$ indicates the time of $R$th failure in a sample of size $n$. Under the Type-II hybrid censoring scheme, we have one of the following three types of observations:\\
Case I: $\{x_{1:n}<\cdots<x_{R:n}\}$ if $x_{R:n}>T.$\\
Case II: $\{x_{1:n}<\cdots<x_{d:n}<T<x_{d+1:n}\}$ if $R\leq d<n$ and $x_{d:n}<T<x_{d+1:n}.$\\
Case III: $\{x_{1:n}<\cdots<x_{n:n}<T\},$ \\
where $x_{1:n}<\cdots<x_{R:n}$ denote the observed ordered failure times of the experimental units. A schematic illustration of the hybrid censoring scheme is presented in Figure \ref{fig1}.

\vspace{1.5cm}
\hspace{1cm}
\small
\begin{picture}(10,30)(10,15)
\put(25,25){\vector(1,0){220}}
\put(35,25){\vector(1,1){20}}
\put(70,25){\vector(1,1){20}}
\put(200,25){\vector(1,1){20}}
\put(25,15){\scriptsize $X_{1:n}$}
\put(35,50){\scriptsize 1-st Failure}
\put(60,15){\scriptsize $X_{2:n}$}
\put(90,50){\scriptsize 2-nd Failure}
\put(190,15){\scriptsize $X_{R:n}$}
\put(180,50){\scriptsize R-th Failure (Experiment Stops)}
\put(165,15){\scriptsize $T$}
\put(122,20){\scriptsize $...$}
\put(300,25){\scriptsize Case I}
\put(35,25){\circle*{1}}
\put(70,25){\circle*{1}}
\put(170,25){\circle*{1}}
\end{picture}

\vspace{1.5cm}
\hspace{1cm}
\small
\begin{picture}(10,30)(10,15)
\put(25,25){\vector(1,0){220}}
\put(35,25){\vector(1,1){20}}
\put(70,25){\vector(1,1){20}}
\put(155,25){\vector(1,1){20}}
\put(200,25){\vector(1,1){20}}
\put(25,15){\scriptsize $X_{1:n}$}
\put(35,50){\scriptsize 1-st Failure}
\put(60,15){\scriptsize $X_{2:n}$}
\put(90,50){\scriptsize 2-nd Failure}
\put(110,20){\scriptsize $...$}
\put(145,15){\scriptsize $X_{d:n}$}
\put(150,50){\scriptsize d-th Failure }
\put(195,15){\scriptsize $T$}
\put(215,50){\scriptsize Experiment Stops}
\put(300,25){\scriptsize Case II}
\put(35,25){\circle*{1}}
\put(70,25){\circle*{1}}
\put(200,25){\circle*{1}}
\end{picture}

\vspace{1.5cm}
\hspace{1cm}
\small
\begin{picture}(10,30)(10,15)
\put(25,25){\vector(1,0){220}}
\put(35,25){\vector(1,1){20}}
\put(70,25){\vector(1,1){20}}
\put(155,25){\vector(1,1){20}}
\put(25,15){\scriptsize $X_{1:n}$}
\put(35,50){\scriptsize 1-st Failure}
\put(60,15){\scriptsize $X_{2:n}$}
\put(90,50){\scriptsize 2-nd Failure}
\put(110,20){\scriptsize $...$}
\put(145,15){\scriptsize $X_{n:n}$}
\put(150,50){\scriptsize n-th Failure (Experiment Stops) }
\put(195,15){\scriptsize $T$}
\put(300,25){\scriptsize Case III}
\put(35,25){\circle*{1}}
\put(70,25){\circle*{1}}
\put(200,25){\circle*{1}}
\put(10,-5){\scriptsize \label{fig1}Figure 1: A schematic presentation for Type-II hybrid censored scheme.}
\end{picture}

\vspace{0.5in}
\par
In this article, we consider the analysis of Type-II hybrid censored lifetime data when the lifetime of each experimental unit follows a two-parameter weighted exponential (WE) distribution. This distribution was originally proposed by Gupta and Kundu \cite{gup1}. The two-parameter WE distribution with the shape and scale parameters $\alpha>0$ and $\lambda>0$, respectively, has the probability density function (pdf) as:
\begin{equation}\label{we}
f_{WE}(x;\alpha,\lambda)=\frac{\alpha+1}{\alpha}\lambda e^{-\lambda x}(1-e^{-\alpha\lambda x});~~~~~~x>0.
\end{equation}
 We denote a two-parameter WE distribution with the pdf (\ref{we}) by $WE(\alpha,\lambda)$ and the corresponding cumulative distribution function (cdf) by $F_{WE}(x;\alpha,\lambda)$.
\par
The aim of this article is two fold. First, we try to earn the MLE's of the unknown parameters. It is observed that the maximum likelihood
estimators can be obtained implicitly by solving two nonlinear equations, but they cannot be obtained in closed form. So MLE's of parameters are derived numerically. Newton-Raphson algorithm is one of the standard methods to determine the MLE's of the parameters. To employ the algorithm, second derivatives of the log-likelihood are required for all iterations. The EM algorithm is a very powerful tool in handling the incomplete data problem see Dempster et al. \cite{dem} and McLachlan and Krishnan \cite{mcl}. Then we use the EM algorithm to compute the MLE's. We also evaluate the observed Fisher information matrix using the missing information principle which have been used to obtained asymptotic confidence intervals of the unknown parameters. The second aim of this article is to provide the Bayes inference for the unknown parameters for Type-II hybrid censored data. It is observed that Bayes estimators can not be obtained explicitly, we provide two approximations namely Lindley's approximation and Gibbs sampling procedure. So we use the Gibbs sampling procedure to compute the Bayes estimators, and the HPD confidence intervals. We compare the performances of the different methods by Monte Carlo simulations, and for illustrative purposes we have analyzed one real data set.
\par
The rest of the article is arranged as follows. In Section 2, we provide The MLE's of the unknown parameters. Fisher information matrix is evaluated in Section 3. Using Lindley's approximation and Gibbs sampling we obtain Bayes estimators and HPD confidence intervals for the parametes in Section 4. Simulation results are presented in section 5. We verify our theoretical results via analyzing data set in Section 6.

\section{Maximum likelihood estimators}
In this section, we study MLEs of the model parameters $\alpha$ and $\lambda$ for $\mbox{WE}(\alpha,\lambda)$ distribution with density function:
$$f(x)=\frac{\alpha+1}{\alpha}\lambda e^{-\lambda x}(1-e^{-\alpha\lambda x}):~\alpha,\lambda,x>0.$$
For simplicity, we apply a re-parametrization as $\alpha$ and $\beta=\alpha \lambda$. By this, the $\mbox{WE}(\alpha,\lambda)$ distribution can be written as:
\begin{equation}
\label{bl}f(x)=\frac{\alpha+1}{\alpha^2}\beta e^{-\frac{\beta}{\alpha}x}(1-e^{-\beta x}):~\alpha,\beta,x>0.
\end{equation}
The likelihood function in Case I is given by
\begin{equation}
\label{a}\hspace{-1.7in}L(\alpha,\lambda)=\frac{n!}{(n-R)!}\Pi_{i=1}^Rf(x_i)(1-F(x_{(R)}))^{(n-R)},
\end{equation}
for Case II,
\begin{equation}
\label{b}\hspace{-1.9in}L(\alpha,\lambda)=\frac{n!}{(n-d)!}\Pi_{i=1}^df(x_i)(1-F(T))^{(n-d)},
\end{equation}
and for case III,
\begin{equation}
\label{III}\hspace{-1.9in}L(\alpha,\lambda)=\Pi_{i=1}^nf(x_i),
\end{equation}
where $f(x)$ is presented by (\ref{bl}), so
$$\hspace{-2.2in}F(x)=1-\frac{1}{\alpha}e^{-\frac{\beta}{\alpha}x}(\alpha+1-e^{-\beta x}).$$
We present likelihood functions (\ref{a}), (\ref{b}) and (\ref{III}) by:
\begin{equation}
\label{c}\hspace{-1.9in}L(\alpha,\lambda)=\frac{n!}{(n-r)!}\Pi_{i=1}^rf(x_i)(1-F(c))^{(n-r)},
\end{equation}
where
\begin{equation}
r=\left\{
\begin{array}{rl}
R & \text{for Case I}\\
d & \text{for Case II}\\
n & \text{for Case III},
\end{array}
\right.
\end{equation}
\mbox{and}
\begin{equation}
c=\left\{
\begin{array}{rl}
x_{R:n} & \text{for Case I}\\
T       & \text{for Cases II and III.}
\end{array}
\right.
\end{equation}
Taking the logarithm of Equation \ref{c}, we obtain
$$l(\alpha,\lambda)=r\ln(\alpha+1)-(n+r)\ln(\alpha)+r\ln(\beta)-\frac{\beta}{\alpha}\sum_{i=1}^rx_i+\sum_{i=1}^r{\ln(1-e^{-\beta x_i})}$$
\begin{equation}
\label{1}+(n-r)(-\frac{\beta}{\alpha})c+(n-r)\ln(\alpha+1-e^{-\beta c}),
\end{equation}
then the normal equations are

\begin{equation}
\displaystyle
\left\{
\begin{array}c
\hspace{-.85 in}\frac{\partial l}{\partial\alpha}=\frac{r}{\alpha+1}-\frac{n+r}{\alpha}+\frac{\beta}{\alpha^2}(\sum_{i=1}^rx_i+(n-r)c)+\frac{(n-r)}{\alpha+1-e^{-\beta c}}\\
\frac{\partial l}{\partial\beta}=\frac{r}{\beta}-\frac{1}{\alpha}(\sum_{i=1}^rx_i+(n-r)c)+\sum_{i=1}^r\frac{x_ie^{-\beta x_i}}{1-e^{-\beta x_i}}+(n-r)\frac{ce^{-\beta c}}{\alpha+1-e^{-\beta c}}
.\end{array}
\right.
\end{equation}

Maximum likelihood estimators can be secured by solving these equations, but they cannot  be expressed explicitly. So we use EM algorithm to compute them. The advantage of this method is that it is convergence for any initial value fast enough.\\


\subsection{EM algorithm}
The EM algorithm, originally proposed by Dempster et al. \cite{dem}, is a very powerful tool for handling
the incomplete data problem.\\
Let us symbolize the observed and
the censored data by $X=(X_{1:n},\cdots,X_{r:n})$ and $Z=(Z_1,\cdots,Z_{n-r})$, respectively. Here
for a given r, $(Z_1,\cdots,Z_{n-r})$ are not observable. The censored data vector $Z$ can
be thought of as missing data. The combination of $W=(X,Z)$ forms the whole
data set.
In next we follow the method Kundu and Pradhan \cite{kun2} for missing data introducing. \\ If we denote the log-likelihood function of the uncensored data set by
$$l_c(\alpha,\beta)=n\ln(\alpha+1)-2n\ln(\alpha)+n\ln(\beta)-\frac{\beta}{\alpha}\left(\sum_{i=1}^rX_i+\sum_{i=1}^{n-r}Z_i\right)$$
\begin{equation}
\label{2}+\left(\sum_{i=1}^r{\ln(1-e^{-\beta x_i})}+\sum_{i=1}^r{\ln(1-e^{-\beta z_i})}\right).
\end{equation}
For the E-step of the EM algorithm, one needs to compute the pseudo log-likelihood function as $l_s(\alpha,\beta)=E(l_c(\alpha,\beta|X)).$
Therefor,
$$l_s(\alpha,\beta)=n\ln(\alpha+1)-2n\ln(\alpha)+n\ln(\beta)-\frac{\beta}{\alpha}\left(\sum_{i=1}^rX_i\right)+
\sum_{i=1}^r{\ln(1-e^{-\beta x_i})}$$
$$-\frac{\beta}{\alpha}(n-r)A(c;\alpha,\beta)+(n-r)B(c;\alpha,\beta),$$
where
$$A[c;\alpha,\beta]=E(Z_i|Z_i>c)
~~\mbox{and}~~
B(c;\alpha,\beta)=E[\ln(1-e^{-\beta Z_i})|Z_i>c],$$
and they are obtained in Appendix A.

Now the M-step includes the maximization of the pseudo log-likelihood
function \ref{2}. Therefore, if at the kth stage, the estimation of $(\alpha,\beta)$ is $(\hat{\alpha}_k,\hat{\beta}_k)$, then $(\hat{\alpha}_{k+1},\hat{\beta}_{k+1})$
can be obtained by maximizing
$$g(\alpha,\beta)=n\ln(\alpha+1)-2n\ln(\alpha)+n\ln(\beta)-\frac{\beta}{\alpha}\left(\sum_{i=1}^rX_i\right)+
\sum_{i=1}^r{\ln(1-e^{-\beta x_i})}$$
\begin{equation}
\label{3}-\frac{\beta}{\alpha}(n-r)A(c;\hat{\alpha}_k,\hat{\beta}_k)+(n-r)B(c;\hat{\alpha}_k,\hat{\beta}_k)
\end{equation}
Note that the maximization of \ref{3} can be earned quite effectively by the similar method proposed by Gupta and Kundu \cite{gup2}. First, $\hat{\beta}_{k+1}$ can be obtain by solving a fixed-point type equation
$$h(\beta)=\beta.$$
The function $h(\beta)$ is defined
$$h(\beta)=n\left[\frac{B}{\hat{\alpha}(\beta)}-\sum_{i=1}^r\frac{x_ie^{-\beta x_i}}{1-e^{-\beta x_i}}\right]^{-1}$$
where
$$B=\sum_{i=1}^rx_i+(n-r)A(c,\hat{\alpha}_k,\hat{\beta}_k)$$ and
$$\hat{\alpha}(\beta)=\frac{\sqrt{(\beta B-2n)^2+4n\beta B}+(\beta B-2n)}{2n}.$$
One can follow iteration method.
Once $\hat{\beta}_{k+1}$ is determined, $\hat{\alpha}_{k+1}$ can be evaluated as $\hat{\alpha}_{k+1}=\hat{\alpha}(\hat{\beta}_{k+1})$.
\par For the estimation of $\lambda$, we can use the invariance property maximum likelihood estimators and obtain $\hat{\lambda}$ as follow:
$$\hat{\lambda}=\frac{\hat{\beta}}{\hat{\alpha}}.$$

\section{Fisher Information matrices}
 One of the advantages of using EM algorithm is that presents a measure of information in censored data through the missing information principle. Louis \cite{lou} improved a procedure for extracting the observed information matrix. In this section, we display the observed Fisher information matrix by using the missing value principles of Louis \cite{lou}. The observed Fisher information matrix can be used to build the asymptotic confidence intervals.
\\
Using the notations: $\theta=(\alpha,\beta)$, X=observed data, W=complete data, $I_X(\theta)$=observed information, $I_W(\theta)$=complete information and
$I_{W|X}(\theta)$=missing information, follow the relation to
\begin{equation}
\label{Ix}I_X(\theta)=I_W(\theta)-I_{W|X}(\theta),
\end{equation}
to evaluate $I_X(\theta).$\\
Complete information and the missing information are given respectively as:
$$I_W(\theta)=-E\left[\frac{\partial^2L_c(W;\theta)}{\partial\theta^2}\right]$$
\mbox{and}
\begin{equation}
\label{IWX}I_{W|X}(\theta)=-(n-r)E\left[\frac{\partial^2\ln f_Z(z|y,\theta)}{\partial\theta^2}\right].
\end{equation}
As the dimension of $\theta$ is 2, $I_X(\theta)$ and $I_{W|X}(\theta)$ are both of the order $2\times2$.\par
The elements of matrix $I_W(\theta)$ for complete data set are presented in Gupta and Kundu \cite{gup1}. They re-parametrized $\mbox{WE}(\alpha,\lambda)$ distribution as $\lambda$ and $\beta=\alpha\lambda$.\\ We report $I_W(\theta)$ which have been evaluated by them here as:
\begin{equation*}
I_W(\theta)=\left[
\begin{array} {cc}
a_{11} & a_{12}\\
a_{21} & a_{22}
\end{array}
\right]
\end{equation*}
where
$$\hspace{-3.3in}\displaystyle a_{11}=\frac{n}{(\beta+\lambda)^2}+\frac{n}{\lambda^2}.$$
$$\hspace{-3.4in}\displaystyle a_{12}=a_{21}=\frac{n}{(\beta+\lambda)^2}.$$
$$\hspace{-2.4in} \displaystyle a_{22}=\frac{n}{(\beta+\lambda)^2}-\frac{n}{\beta^2}+\frac{n\lambda(\beta+\lambda)A}{\beta^4},$$
in which $\displaystyle A=\int_0^1\frac{(\ln (1-y))^2(1-y)^{\frac{\lambda}{\beta}}}{y}dy.$\\
On the other hand, with the above re-parametrization and by using (\ref{IWX}), one can easily verify
\begin{equation*}
I_{W|X}(\theta)=(n-r)\left[
\begin{array} {cc}
b_{11}(c;\alpha,\beta) & b_{12}(c;\alpha,\beta)\\
b_{21}(c;\alpha,\beta) & b_{22}(c;\alpha,\beta)
\end{array}
\right],
\end{equation*}
where
$$\hspace{-1.in}\displaystyle b_{11}(c;\alpha,\beta)=\frac{1}{(\beta+\lambda)^2}+\frac{2\beta}{\lambda^3(\frac{\beta}{\lambda}+1-e^{-\beta c})}-\frac{\beta^2}{\lambda^4(\frac{\beta}{\lambda}+1-e^{-\beta c})^2}.$$
$$\hspace{-0.2in}\displaystyle b_{12}(c;\alpha,\beta)=b_{21}(c;\alpha,\beta)=\frac{1}{(\beta+\lambda)^2}-\frac{1}{\lambda^2(\frac{\beta}{\lambda}+1-e^{-\beta c})}+\frac{\beta(\frac{1}{\lambda}+ce^{-\beta c})}{\lambda^2(\frac{\beta}{\lambda}+1-e^{-\beta c})^2}.$$
$$\hspace{0in}\displaystyle b_{22}(c;\alpha,\beta)=\frac{1}{(\beta+\lambda)^2}-\frac{c^2e^{-\beta c}}{\frac{\beta}{\lambda}+1-e^{-\beta c}}-\frac{(\frac{1}{\lambda}+c-e^{-\beta c})^2}{(\frac{\beta}{\lambda}+1-e^{-\beta c})^2}+\frac{(\beta+\lambda)B}{\beta^3e^{-\lambda c}(\frac{\beta}{\lambda}+1-e^{-\beta c})},$$
in which $B=\int_{1-e^{-\beta c}}^1\frac{(\ln (1-y))^2(1-y)^{\frac{\lambda}{\beta}}}{y}dy.$
\par Now, $I_X(\theta)$ can be computed by (\ref{Ix}). The asymptotic variance-covariance matrix of $\hat{\theta}$ can be obtained by inverting $I_X(\theta)$. We use this matrix to secure the asymptotic confidence intervals for $\lambda$ and $\beta$. To obtain the asymptotic confidence interval for $\alpha$, we use the non-parametric bootstrap method \cite{zog}.

\section{Bayes Estimators and Confidence Intervals}
In this section, we study Bayes estimators for parameters $\alpha$ and $\lambda$ under symmetric loss functions. A very well known symmetric loss function is the squared error which is defined as:
$L(f(\mu),\hat{f}(\mu))=(\hat{f}(\mu)-f(\mu))^2,$ with $\hat{f}(\mu)$ being an estimate of $f(\mu)$. Here $f(\mu)$ denotes
some function of $\mu$. Bayes estimators, say $\hat{f}_{Bayes}(\mu)$, is evaluated by the
posterior mean of $f(\mu)$.\par
Let $x=(x_{1:n},\cdots,x_{r:n})$ be an observed sample from the hybrid censoring scheme, drawn from a $\mbox{WE}(\alpha,\lambda)$ distribution. We apply re-parametrization as $\alpha$ and $\beta=\alpha\lambda$. So the likelihood function becomes
$$L(\alpha,\beta|\b{x})\varpropto\frac{(\alpha+1)^r}{\alpha^{n+r}}\beta^re^{-\frac{\beta}{\alpha}\left(\sum_{i=1}^rx_i+c(n-r)\right)}(\alpha+1-e^{-\beta c})^{n-r}\Pi_{i=1}^r(1-e^{-\beta x_i}),$$
and $\log$-likelihood function:
\begin{equation*}
l(\alpha,\beta|\b{x})=r\log(\alpha+1)-(n+r)\log\alpha+r\log(\beta)-\frac{\beta}{\alpha}\left(\sum_{i=1}^rx_i+c(n-r)\right)
\end{equation*}
\begin{equation}\label{log}
+(n-r)\log(\alpha+1-e^{-\beta c})+\sum_{i=1}^r\log(1-e^{-\beta x_i}).
\end{equation}
It is assumed that $\beta$ and $\alpha$ have the following independent
gamma priors:
$$\pi_1(\beta)\varpropto\beta^{w_2-1}e^{-\beta w_1},~~~~~~\beta>0,$$
$$\pi_2(\alpha)\varpropto\alpha^{w_4-1}e^{-\alpha w_3},~~~~~~\alpha>0.$$
So, the joint prior distribution of $\alpha$ and $\beta$ is of the form
$$\pi(\alpha,\beta)\varpropto\alpha^{w_4-1}e^{-\alpha w_3}\beta^{w_2-1}e^{-\beta w_1},~~~~~~~\alpha>0,~\beta>0,~w_1>0,~w_2>0,~w_3>0,~w_4>0.$$
Then the posterior  distribution $\alpha$ and $\beta$ can be written as
$$\pi(\alpha,\beta|\b{x})=\frac{1}{k}\alpha^{w_4-n-r-1}\beta^{w_2+r-1}(\alpha+1)^re^{-\alpha w_3}e^{-\beta w_1}e^{-\frac{\beta}{\alpha}\left(\sum_{i=1}^rx_i+c(n-r)\right)}$$
\begin{equation}
\times(\alpha+1-e^{-\beta c})^{n-r}\Pi_{i=1}^r(1-e^{-\beta x_i})\label{pos}
\end{equation}
where
$$k=\int_0^\infty\int_0^\infty\alpha^{w_4-n-r-1}\beta^{w_2+r-1}(\alpha+1)^re^{-\alpha w_3}e^{-\beta w_1}e^{-\frac{\beta}{\alpha}\left(\sum_{i=1}^rx_i+c(n-r)\right)}$$
$$\times(\alpha+1-e^{-\beta c})^{n-r}\Pi_{i=1}^r(1-e^{-\beta x_i})d\alpha d\beta.$$
Now the Bayes estimators of $\alpha$ and $\beta$ under the squared error loss function L are respectively obtained as:
$$\hat{\alpha}_{Bayes}=E[\alpha|\b{x}]=\frac{1}{k}\int_0^\infty\int_0^\infty\alpha^{w_4-n-r}\beta^{w_2+r-1}(\alpha+1)^re^{-\alpha w_3}e^{-\beta w_1}$$$$\times e^{-\frac{\beta}{\alpha}\left(\sum_{i=1}^rx_i+c(n-r)\right)}(\alpha+1-e^{-\beta c})^{n-r}\Pi_{i=1}^r(1-e^{-\beta x_i})d\alpha d\beta,$$
and
$$\hat{\beta}_{Bayes}=E[\beta|\b{x}]=\frac{1}{k}\int_0^\infty\int_0^\infty\alpha^{w_4-n-r-1}\beta^{w_2+r}(\alpha+1)^re^{-\alpha w_3}e^{-\beta w_1}$$$$\times e^{-\frac{\beta}{\alpha}\left(\sum_{i=1}^rx_i+c(n-r)\right)}(\alpha+1-e^{-\beta c})^{n-r}\Pi_{i=1}^r(1-e^{-\beta x_i})d\alpha d\beta.$$
Since $\lambda$ is a function of $\alpha$ and $\beta$, then one can obtain the posterior density function of $\lambda$ and so the Bayes estimator of $\lambda$ under the squared error loss function $L$ as:
$$\hat{\lambda}_{Bayes}=E[\lambda|\b{x}]=\frac{1}{k}\int_0^\infty\int_0^\infty u^{w_4+w_2-n-1}\lambda^{w_2+r-1}(1+u)^{r}e^{-u(w_3+\lambda w_1)}$$$$
\times e^{-\lambda\left(\sum_{i=1}^rx_i+c(n-r)\right)}(u+1-e^{-\lambda u c})^{n-r}\Pi_{i=1}^r(1-e^{-\lambda u x_i})du d\lambda.$$
As these estimators can not be evaluated explicitly, so we adopt two different procedures to approximate them:
\begin{itemize}
\item{ Lindley approximation,}
\item{ MCMC method}.
\end{itemize}
\subsection{Lindley approximation method}
\par
In previous section, based on Type-II hybrid censored scheme we obtained the Bayes estimators of $\alpha$, $\beta$ and $\lambda$ against squared error loss function $L$. It is easily observed that theses estimators have not explicit closed forms. For these evaluation, numerical techniques are required. One of the most numerical techniques is Lindley's method (see \cite{lin}), that for these estimators can be describe as follows.
In general, Bayes estimator of $u(\alpha,\beta)$ as a function of $\alpha$ and $\beta$ is identified:
$$I(\b{x})=\frac{\int_0^\infty\int_0^\infty u(\alpha,\beta)e^{l(\alpha,\beta|\b{x})+\rho(\alpha,\beta)}d\alpha d\beta}{\int_0^\infty\int_0^\infty e^{l(\alpha,\beta|\b{x})+\rho(\alpha,\beta)}d\alpha d\beta},$$
where $l(\alpha,\beta|\b{x})$ is $\log$-likelihood function (defined by \ref{log}) and $\rho(\alpha,\beta)=\log\pi(\alpha,\beta)$.
\\ By the Lindley's method $I(\b{x})$ can be approximated as:
$$I(\b{x})=u(\hat{\alpha},\hat{\beta})+\frac{1}{2}[(\hat{u}_{\alpha\alpha}+2\hat{u}_{\alpha}\hat{\rho}_{\alpha})\hat{\sigma}_{\alpha\alpha}+
(\hat{u}_{\beta\alpha}+2\hat{u}_{\beta}\hat{\rho}_{\alpha})\hat{\sigma}_{\beta\alpha}+
(\hat{u}_{\alpha\beta}+2\hat{u}_{\alpha}\hat{\rho}_{\beta})\hat{\sigma}_{\alpha\beta}$$$$
+(\hat{u}_{\beta\beta}+2\hat{u}_{\beta}\hat{\rho}_{\beta})\hat{\sigma}_{\beta\beta}]+\frac{1}{2}[(\hat{u}_{\alpha}\hat{\sigma}_{\alpha\alpha}+
\hat{u}_{\beta}\hat{\sigma}_{\alpha\beta})(\hat{l}_{\alpha\alpha\alpha}\hat{\sigma}_{\alpha\alpha}+
\hat{l}_{\alpha\beta\alpha}\hat{\sigma}_{\alpha\beta}+\hat{l}_{\beta\alpha\alpha}\hat{\sigma}_{\beta\alpha}$$$$
+\hat{l}_{\beta\beta\alpha}\hat{\sigma}_{\beta\beta})+(\hat{u}_{\alpha}\hat{\sigma}_{\beta\alpha}+
\hat{u}_{\beta}\hat{\sigma}_{\beta\beta})(\hat{l}_{\beta\alpha\alpha}\hat{\sigma}_{\alpha\alpha}
+\hat{l}_{\alpha\beta\beta}\hat{\sigma}_{\alpha\beta}+\hat{l}_{\beta\alpha\beta}\hat{\sigma}_{\beta\alpha}
+\hat{l}_{\beta\beta\beta}\hat{\sigma}_{\beta\beta})],$$
where $\hat{\alpha}$ and $\hat{\beta}$ are the MLE's of $\alpha$ and $\beta$ respectively. Also, $u_{\alpha\alpha}$ is the second derivative of the function $u(\alpha,\beta)$ with the respect to $\alpha$ and $\hat{u}_{\alpha\alpha}$ valued of $u_{\alpha\alpha}$ at $(\hat{\alpha},\hat{\beta}).$ Other expressions can be calculated with following definitions:
$$\hspace{-1in}\hat{l}_{\alpha\alpha}=\frac{\partial^2l}{\partial\alpha^2}\left|_{\alpha=\hat{\alpha},\beta=\hat{\beta}}=\frac{-r}{(\hat{\alpha}+1)^2}+\frac{n+r}{\hat{\alpha}^2}-\frac{2\hat{\beta} A}{\hat{\alpha}^3}-\frac{n-r}{(\hat{\alpha}+1-e^{-\hat{\beta} c})^2}\right.,$$
$$\hspace{-2in}\hat{l}_{\alpha\beta}=\frac{\partial^2l}{\partial\alpha\partial\beta}\left|_{\alpha=\hat{\alpha},\beta=\hat{\beta}}=\hat{l}_{\beta\alpha}=\frac{A}{\hat{\alpha}^2}-\frac{c(n-r)e^{-\hat{\beta} c}}{(\hat{\alpha}+1-e^{-\hat{\beta} c})^2}\right.,$$
$$\hspace{-.2in}\hat{l}_{\alpha\beta\alpha}=\frac{\partial^3l}{\partial\alpha\partial\beta\partial\alpha}\left|_{\alpha=\hat{\alpha},\beta=\hat{\beta}}=\hat{l}_{\beta\alpha\alpha}=\frac{\partial^3l}{\partial\beta\partial\alpha^2}\left|_{\alpha=\hat{\alpha},\beta=\hat{\beta}}=\frac{-2A}{\hat{\alpha}^3}+\frac{2c(n-r)e^{-\hat{\beta} c}}{(\hat{\alpha}+1-e^{-\hat{\beta} c})^3}\right.\right.,$$
$$\hspace{-0.8in}\hat{l}_{\alpha\alpha\alpha}=\frac{\partial^3l}{\partial\alpha^3}\left|_{\alpha=\hat{\alpha},\beta=\hat{\beta}}=\frac{2r}{(\hat{\alpha}+1)^3}-\frac{2(n+r)}{\hat{\alpha}^3}+\frac{6\hat{\beta} A}{\hat{\alpha}^4}+\frac{2(n-r)}{(\hat{\alpha}+1-e^{-\hat{\beta} c})^3}\right.,$$
where $A=\sum_{i=1}^rx_i+(n-r)c,$
$$\hspace{0in}\hat{l}_{\alpha\beta\beta}=\frac{\partial^3l}{\partial\alpha\partial\beta^2}\left|_{\alpha=\hat{\alpha},\beta=\hat{\beta}}=\hat{l}_{\beta\alpha\beta}=\frac{\partial^3l}{\partial\beta\partial\alpha\partial\beta}\left|_{\alpha=\hat{\alpha},\beta=\hat{\beta}}=\frac{c^2(n-r)e^{-\hat{\beta} c}(\hat{\alpha}+1+e^{-\hat{\beta}c})}{(\hat{\alpha}+1-e^{-\hat{\beta} c})^3}\right.\right.,$$
$$\hspace{-.9in}\hat{l}_{\beta\beta}=\frac{\partial^2l}{\partial\beta^2}\left|_{\alpha=\hat{\alpha},\beta=\hat{\beta}}=\frac{-r}{\hat{\beta}^2}-\frac{c^2(n-r)(\hat{\alpha}+1)e^{-\hat{\beta} c}}{(\hat{\alpha}+1-e^{-\hat{\beta} c})^2}-\sum_{i=1}^r\frac{x_i^2e^{-\hat{\beta} x_i}}{(1-e^{-\hat{\beta} x_i})^2}\right.,$$
$$\hspace{-1.9in}\hat{l}_{\beta\beta\alpha}=\frac{\partial^3l}{\partial\beta^2\partial\alpha}\left|_{\alpha=\hat{\alpha},\beta=\hat{\beta}}=\frac{c^2(n-r)e^{-\hat{\beta} c}(\hat{\alpha}+1+e^{-\hat{\beta}c})}{(\hat{\alpha}+1-e^{-\hat{\beta} c})^3}\right.,$$
$$\hspace{0.05in}\hat{l}_{\beta\beta\beta}=\frac{\partial^3l}{\partial\beta^3}\left|_{\alpha=\hat{\alpha},\beta=\hat{\beta}}=\frac{2r}{\hat{\beta}^3}+\frac{c^3(n-r)(\hat{\alpha}+1)e^{-\hat{\beta} c}(\hat{\alpha}+1+e^{-\hat{\beta}c})}{(\hat{\alpha}+1-e^{-\hat{\beta} c})^3}+\sum_{i=1}^r\frac{x_i^3e^{-\hat{\beta} x_i}(1+e^{-\hat{\beta}x_i})}{(1-e^{-\hat{\beta} x_i})^3}\right.,$$
$$\hspace{-3.6in}\hat{\rho}_{\beta}=\frac{\partial\rho}{\partial\beta}\left|_{\beta=\hat{\beta}}=\frac{w_2-1}{\hat{\beta}}-w_1\right.,$$
$$\hspace{-3.6in}\hat{\rho}_{\alpha}=\frac{\partial\rho}{\partial\alpha}\left|_{\alpha=\hat{\alpha}}=\frac{w_4-1}{\hat{\alpha}}-w_3\right.,$$
and we have:
\[ \left( \begin{array}{cc}
\hat{\sigma}_{\alpha\alpha} & \hat{\sigma}_{\alpha\beta} \\
\hat{\sigma}_{\beta\alpha} & \hat{\sigma}_{\beta\beta}  \end{array} \right) = \left( \begin{array}{cc}
-\hat{l}_{\alpha\alpha} & -\hat{l}_{\alpha\beta} \\
-\hat{l}_{\beta\alpha} & -\hat{l}_{\beta\beta} \end{array} \right)^{-1}.\]
With the above defined expressions, we obtain the approximation Bayes estimators.\\
Also we have:
$$u(\alpha,\beta)=\alpha,~~~~u_{\alpha}=1,~~~~u_{\alpha\alpha}=u_{\beta}=u_{\beta\beta}=u_{\alpha\beta}=u_{\beta\alpha}=0,$$
the Bayes estimator of $\alpha$ under the squared error loss function $L$ becomes
$$\hat{\alpha}_{Bayes}=\hat{\alpha}+\frac{1}{2}[2\hat{\rho}_{\alpha}\hat{\sigma}_{\alpha\alpha}+
2\hat{\rho}_{\beta}\hat{\sigma}_{\alpha\beta}+\hat{\sigma}^2_{\alpha\alpha}\hat{l}_{\alpha\alpha\alpha}+
3\hat{\sigma}_{\alpha\alpha}\hat{\sigma}_{\alpha\beta}\hat{l}_{\beta\alpha\alpha}+2\hat{\sigma}^2_{\alpha\beta}\hat{l}_{\alpha\beta\beta}$$$$
+\hat{\sigma}_{\alpha\alpha}\hat{\sigma}_{\beta\beta}\hat{l}_{\beta\beta\alpha}+
\hat{\sigma}_{\alpha\beta}\hat{\sigma}_{\beta\beta}\hat{l}_{\beta\beta\beta}].$$
Proceeding similarly, the Bayes estimator of $\beta$ under $L$ is given by
$$(u(\alpha,\beta)=\beta,~~~~u_{\beta}=1,~~~~u_{\alpha\alpha}=u_{\alpha}=u_{\beta\beta}=u_{\alpha\beta}=u_{\beta\alpha}=0),$$
$$\hat{\beta}_{Bayes}=\hat{\beta}+\frac{1}{2}[2\hat{\rho}_{\alpha}\hat{\sigma}_{\alpha\beta}+
2\hat{\rho}_{\beta}\hat{\sigma}_{\beta\beta}+2\hat{\sigma}^2_{\alpha\beta}\hat{l}_{\beta\alpha\alpha}+
2\hat{\sigma}_{\alpha\beta}\hat{\sigma}_{\beta\beta}\hat{l}_{\alpha\beta\beta}$$
$$+\hat{\sigma}^2_{\beta\beta}\hat{l}_{\beta\beta\beta}+
\hat{\sigma}_{\alpha\alpha}\hat{\sigma}_{\alpha\beta}\hat{l}_{\alpha\alpha\alpha}+
\hat{\sigma}_{\beta\beta}\hat{\sigma}_{\alpha\alpha}\hat{l}_{\beta\alpha\alpha}+
\hat{\sigma}_{\alpha\beta}\hat{\sigma}_{\beta\beta}\hat{l}_{\beta\beta\alpha}].$$
Finally the Bayes estimator of $\lambda$ under $L$ is given by
$$(u(\alpha,\beta)=\frac{\beta}{\alpha},~~u_{\beta}=\frac{1}{\alpha},~~u_{\alpha}=\frac{-\beta}{\alpha^2},~~u_{\alpha\alpha}=\frac{2\beta}{\alpha^3},
~~u_{\beta\beta}=0,~~u_{\alpha\beta}=u_{\beta\alpha}=\frac{-1}{\alpha^2}),$$
$$\hat{\lambda}_{Bayes}=\frac{\hat{\beta}}{\hat{\alpha}}+\frac{1}{2}[(\hat{u}_{\alpha\alpha}+2\hat{u}_{\alpha}\hat{\rho}_{\alpha})\hat{\sigma}_{\alpha\alpha}+
(\hat{u}_{\beta\alpha}+2\hat{u}_{\beta}\hat{\rho}_{\alpha})\hat{\sigma}_{\beta\alpha}+
(\hat{u}_{\alpha\beta}+2\hat{u}_{\alpha}\hat{\rho}_{\beta})\hat{\sigma}_{\alpha\beta}$$$$
+2\hat{u}_{\beta}\hat{\rho}_{\beta}\hat{\sigma}_{\beta\beta}]+\frac{1}{2}[(\hat{u}_{\alpha}\hat{\sigma}_{\alpha\alpha}+
\hat{u}_{\beta}\hat{\sigma}_{\alpha\beta})(\hat{l}_{\alpha\alpha\alpha}\hat{\sigma}_{\alpha\alpha}+
\hat{l}_{\alpha\beta\alpha}\hat{\sigma}_{\alpha\beta}+\hat{l}_{\beta\alpha\alpha}\hat{\sigma}_{\beta\alpha}$$$$
\hspace{0.5in}+\hat{l}_{\beta\beta\alpha}\hat{\sigma}_{\beta\beta})+(\hat{u}_{\alpha}\hat{\sigma}_{\beta\alpha}+
\hat{u}_{\beta}\hat{\sigma}_{\beta\beta})(\hat{l}_{\beta\alpha\alpha}\hat{\sigma}_{\alpha\alpha}
+\hat{l}_{\alpha\beta\beta}\hat{\sigma}_{\alpha\beta}+\hat{l}_{\beta\alpha\beta}\hat{\sigma}_{\beta\alpha}
+\hat{l}_{\beta\beta\beta}\hat{\sigma}_{\beta\beta})].$$
\par
The approximate Bayes estimators of $\alpha$, $\beta$ and $\lambda$ can be obtained using Lindley approximation, but it is not possible to construct highest posterior density (HPD) confidence intervals using this method. Therefore, we suggest the following Markov Chain Monte Carlo (MCMC) method to generate samples from the posterior density function, and in turn to obtain the Bayes estimators, and HPD confidence intervals.
\subsection{Gibbs sampling}
\par
\par
 Here we study the Gibbs sampling method to draw samples from the posterior density function and then compute the Bayes estimators and HPD confidence intervals of $\alpha$, $\beta$ and $\lambda$ under the squared errors loss function.
\par
Let $\b{x}=(x_{1:n},\cdots,x_{r:n})$ be an observed sample from the hybrid censoring scheme, drawn from a $\mbox{WE}(\alpha,\lambda)$ distribution.
From (\ref{pos}), we can write the joint posterior density function of $\alpha$ and $\beta$ given $\b{x}$ as:
$$\pi(\alpha,\beta|\b{x})\propto\alpha^{w_4-n-r-1}\beta^{w_2+r-1}(\alpha+1)^re^{-\alpha w_3}e^{-\beta w_1}\left\{e^{-\frac{\beta}{\alpha}\left(\sum_{i=1}^rx_i+c(n-r)\right)}\right.$$
\begin{equation}
\left.\times(\alpha+1-e^{-\beta c})^{n-r}\Pi_{i=1}^r(1-e^{-\beta x_i})\right\}\label{pos2},
\end{equation}
by this, the posterior density function of $\beta$ given $\alpha$ and $\b{x}$ is
$$\pi(\beta|\alpha,\b{x})\propto\beta^{w_2+r-1}e^{-\beta \left(w_1+\frac{\beta}{\alpha}\sum_{i=1}^rx_i+c(n-r)\right)}
(\alpha+1-e^{-\beta c})^{n-r}\Pi_{i=1}^r(1-e^{-\beta x_i}).$$
{\label{the1}\theorem The conditional distribution of $\beta$ given $\alpha$ and $\b{x}$ is log-concave.}
{\proof See Appendix, part B.}\\
By (\ref{pos2}), the posterior density function of $\alpha$ given $\beta$ and $\b{x}$ is
\begin{equation}
\pi(\alpha|\beta,\b{x})\propto\alpha^{w_4-n-r-1}(\alpha+1)^re^{-\alpha w_3}e^{-\frac{\beta}{\alpha}\left(\sum_{i=1}^rx_i+c(n-r)\right)}(\alpha+1-e^{-\beta c})^{n-r}.\label{albe}
\end{equation}
{\label{th2}\theorem The conditional distribution of $\alpha$ given $\beta$ and $\b{x}$ has a finite maximum point.}
{\proof See Appendix, part C.}
{\corollary With the help of the acceptance rejection principle (see Devroye \cite{dev} for details) and the previous theorem, the generation from (\ref{albe}) can be performed using the WE generator.}\\ \par
Now we use theorems \ref{the1} and \ref{th2} and pursue the idea of Geman and Geman \cite{gem}, and suggest the following scheme.
\\
\begin{itemize}
{\item Step 1) Take some initial value of $\alpha$ and $\beta$, such as $\alpha_0$ and $\beta_0$.}
{\item Step 2) Generate $\alpha_{i+1}$ and $\beta_{i+1}$ from $\pi(\alpha|\beta_i,\b{x})$ and $\pi(\beta|\alpha_i,\b{x})$.}
{\item Step 3) Repeat Step 2, $N$ times.}
{\item Step 4) Obtain Bayes estimators of $\alpha$ and $\beta$ with respect to  a squared error loss function:
$$\hat{\alpha}_{Bayes}=\frac{1}{N-M_1}\sum_{i=M_1+1}^{N}{\alpha_i}~~~~\mbox{and}
~~~~\hat{\beta}_{Bayes}=\frac{1}{N-M_2}\sum_{i=M_2+1}^{N}{\beta_i}$$
 where $M_1$ and $M_2$ are the burn-in periods in generating of $\alpha_i$ and $\beta_i$ respectively.}
{\item Step 5) Obtain the HPD confidence interval of $\alpha$: Order $\alpha_1,\cdots,\alpha_{M_1}$ as $\alpha_{(1)}<\cdots<\alpha_{(M_1)}$ and construct all the $100(1-\eta)\%$ confidence intervals of $\alpha$, as:
$$(\alpha_{(1)},\alpha_{([M_1(1-\eta)])}),\cdots,(\alpha_{([M_1\eta])},\alpha_{(M_1)}),$$
where $[M]$  symbolizes the largest integer less than or equal to $M$. The HPD confidence interval of $\alpha$ is the shortest length interval. Similarly, we can construct a $100(1-\eta)\%$ HPD confidence interval of $\beta$.}
\end{itemize}\par
Finally, using the idea of Chen and Shao \cite{chen}, we can compute the estimation and HPD confidence interval for $\lambda$.

\section{Numerical Experiments}
In this section, we carry out a simulation study to compare the performance of MLE's and Bayes estimators. In all the cases $\alpha=2.5$ and $\lambda=3$ are taken. We estimate the unknown parameters using the MLE, Bayes estimators obtained by Lindley's approximations and also Bayes estimators obtained by using MCMC technique. We compare the performances of different estimators with MSE. We also obtain the average length of the asymptotic confidence intervals and the HPD confidence intervals. \par
For computing the Bayes estimators, it is assumed that $\beta$ and $\alpha$ have $\mbox{Gamma}(w_2,w_1)$ and $\mbox{Gamma}(w_4,w_3)$ priors, respectively. Moreover we use the non-informative priors of both $\beta$ and $\alpha$, by considering $w_1=w_2=w_3=w_4=0$. The Bayes estimators are computed under the squared error loss function and with respect to the above non-informative priors. \par
The simulation is performed for different choices of $n,~R,~T$ values. We replicate the procedure for 1000 times and report the average estimators, the MSE's, the average asymptotic confidence intervals length and the average HPD confidence intervals length from the MCMC technique. The results are reported in Table 1-4. The first and second rows are parameter estimators of $\lambda$ and $\alpha$, respectively.

\begin{table}
\caption{{\small Average estimators, corresponding MSE and average confidence (asymptotic or HPD for Gibbs) length for N=40 , T=1}}
\hspace{0.2in}
{\begin{tabular}{|l|l|l|l|}
\hline
~ & R=25 & R=30 & R=35\\
\hline
{\small\mbox{MLE}} & 2.978(4.623)3.013 & 3.014(1.836)2.993 & 2.909(0.008)2.999 \\
~ & 2.575(0.029)0.945 & 2.574(0.028)0.985 & 2.578(0.027)0.814 \\
\hline
{\small\mbox{Bayes(Lindley)}} & 2.863(0.019) & 2.968(0.019) & 2.863(0.019) \\
~ & 2.573(0.005) & 2.570(0.005) & 2.574(0.005) \\
\hline
{\small\mbox{Bayes(Gibbs)}} & 2.987(0.246)1.829 & 2.979(0.179)1.516 & 2.944(0.159)1.451 \\
~ & 2.472(0.055)0.757 & 2.482(0.053)0.751 & 2.489(0.053)0.759 \\
\hline
\end{tabular}}
\end{table}

\begin{table}
\caption{{\small Average estimators, corresponding MSE and average confidence (asymptotic or HPD for Gibbs) length for N=40 , T=2}}
\hspace{0.2in}
{\begin{tabular}{|l|l|l|l|}
\hline
~ & R=25 & R=30 & R=35\\
\hline
{\small\mbox{MLE}} & 3.192(0.037)2.865 & 3.145(0.021)2.899 & 2.908(0.008)2.889 \\
~ & 2.624(0.046)0.972 & 2.591(0.033)0.891 & 2.585(0.032)0.894 \\
\hline
{\small\mbox{Bayes(Lindley)}} & 2.968(0.001) & 2.968(0.001) & 2.963(0.001) \\
~ & 2.621(0.015) & 2.576(0.006) & 2.581(0.006) \\
\hline
{\small\mbox{Bayes(Gibbs)}} & 3.148(0.216)1.544 & 3.004(0.166)1.217 & 2.956(0.146)1.354 \\
~ & 2.479(0.054)0.0752 & 2.482(0.054)0.748 & 2.481(0.053)0.758 \\
\hline
\end{tabular}}
\end{table}

\begin{table}
\caption{{\small Average estimators, corresponding MSE and average confidence (asymptotic or HPD for Gibbs) length for N=50 , T=1}}
\hspace{0.2in}
{\begin{tabular}{|l|l|l|l|}
\hline
~ & R=35 & R=40 & R=45\\
\hline
{\small\mbox{MLE}} & 3.287(0.083)2.669 & 3.079(0.006)2.667 & 3.041(0.002)2.735 \\
~ & 2.556(0.019)0.731 & 2.555(0.018)0.763 & 2.463(0.018)0.784 \\
\hline
{\small\mbox{Bayes(Lindley)}} & 2.979(0.000) & 2.984(0.000) & 2.981(0.000) \\
~ & 2.554(0.003) & 2.552(0.003) & 2.561(0.003) \\
\hline
{\small\mbox{Bayes(Gibbs)}} & 2.958(0.269)1.661 & 3.066(0.247)1.643 & 2.966(0.234)1.318 \\
~ & 2.484(0.054)0.759 & 2.485(0.054)0.757 & 2.485(0.052)0.757 \\
\hline
\end{tabular}}
\end{table}

\begin{table}
\vspace{-.05in}
\caption{{\small Average estimators, corresponding MSE and average confidence (asymptotic or HPD for Gibbs) length for N=50 , T=2}}
\hspace{0.2in}
{\begin{tabular}{|l|l|l|l|}
\hline
~ & R=35 & R=40 & R=45\\
\hline
{\small\mbox{MLE}} & 3.282(0.079)2.587 & 2.921(0.006)2.579 & 2.930(0.005)2.577 \\
~ & 2.557(0.021)0.751 & 2.558(0.022)0.745 & 2.556(0.018)0.783 \\
\hline
{\small\mbox{Bayes(Lindley)}} & 2.973(0.000) & 2.972(0.000) & 2.964(0.001) \\
~ & 2.556(0.003) & 2.555(0.003) & 2.553(0.002) \\
\hline
{\small\mbox{Bayes(Gibbs)}} & 3.194(0.214)1.418 & 3.083(0.155)1.349 & 3.001(0.045)0.504 \\
~ & 2.483(0.053)0.754 & 2.495(0.053)0.754 & 2.493(0.054)0.759 \\
\hline
\end{tabular}}
\end{table}
From Tables 1-4, it is observed that for fixed N and T as R increases, the MSE decrease. The performances of the MLE's and Bayes estimators are very similar in all aspects. The average HPD confidence lengths are smaller than the average asymptotic lengths in all the cases considered. Finally it should be mentioned that Bayes estimators are most computationally expensive followed by MLE's.
\section{Data Analysis}
In this section, we demonstrate one data set for illustrative purposes. It has been studied by Gupta and Kundu \cite{gup1} that the $\mbox{WE}(\alpha,\lambda)$ distribution can be used quite to analyze them and MLE's of $\alpha$ and $\lambda$ are 1.6232 and 0.0138 respectively. The data set was studied by Bjerkedal \cite{bje} and is given below:\\
12 15 22 24 24 32 32 33 34 38 38 43 44 48 52 53 54 54 55 56 57 58 58 59 60 60 60 60 61 62 63
65 65 67 68 70 70 72 73 75 76 76 81 83 84 85 87 91 95 96 98 99 109 110 121 127 129 131 143
146 146 175 175 211 233 258 258 263 297 341 341 376.
\\
We use them and create the following two sampling schemes:
$$\mbox{Scheme 1:}~ T=300,~R=60,$$
$$\mbox{Scheme 2:}~ T=250,~R=65,$$
Now for scheme 1, MLE of $\beta,~\alpha~\mbox{and}~\lambda$ are $(0.0239,1.7715,0.0135)$ and Bayes estimators with assumed non-informative priors, i.e., $w_1=w_2=w_3=w_4=0$ with Lindley approximation and Gibbs sampling method are $(0.0198,1.5019,0.0147)$ and $(0.0256,2.0372,0.0138)$ respectively. The $95\%$ confidence intervals based on MLE and Bayes estimators of $\beta,~\alpha~\mbox{and}~\lambda$ are $$\{(0,0.0664),(0.0073,4.3224),(0.0075,0.0195)\}$$ and $$\{(0.0161,0.0350),(1.1000,2.9996),(0.0053,0.0256)\}$$ respectively. \\
For scheme 2, MLE of $\beta,~\alpha~\mbox{and}~\lambda$ are $(0.0255,1.9390,0.0132)$ Bayes estimators with assumed non-informative priors, i.e., $w_1=w_2=w_3=w_4=0$ with Lindley approximation and Gibbs sampling method are $(0.0223,1.7395,0.0142)$ and $(0.0254,2.0835,0.0131)$ respectively. The $95\%$ confidence intervals based on MLE and Bayes estimators of $\beta,~\alpha~\mbox{and}~\lambda$ are $$\{(0,0.0677),(0.0024,4.5216),(0.0075,0.0189))\}$$ and $$\{(0.0170,0.0340),(1.2242,2.9353),(0.0085,0.0222)\}$$ respectively.\par
Because we see the effect of the hyper parameters on the Bayes estimators and also on confidence intervals, we take the following informative priors $w_1=3,~w_2=1.5,~w_3=0.01,~w_4=1.$ \\
Based on this, for scheme 1, Bayes estimators of $\beta,~\alpha~\mbox{and}~\lambda$ with Lindley approximation and Gibbs sampling method are $(0.0233,1.8233,0.0142)$ and $(0.0255,1.8205,0.0150)$ respectively. The $95\%$ confidence intervals based on Bayes estimators of $\beta,~\alpha~\mbox{and}~\lambda$ are $$\{(0.0162,0.0350),(1.0042,2.6231),(0.0062,0.0269)\}.$$
For scheme 2, Bayes estimators of $\beta,~\alpha~\mbox{and}~\lambda$ with Lindley approximation and Gibbs sampling method are $(0.0231,1.8162,0.0141)$ and $(0.0246,2.0799,0.0127)$ respectively. The $95\%$ confidence intervals based on Bayes estimators of $\beta,~\alpha~\mbox{and}~\lambda$ are $$\{(0.0170,0.0322),(1.2043,2.9135),(0.0060,0.0217)\}.$$
 \par
 We plot all the different estimated density functions with non-informative priors and informative priors in Figure 1 and Figure 2. \label{fig3}\input{epsf}
\epsfxsize=3in \epsfysize=2in
 \begin{figure}
\centerline{\epsfxsize=5in \epsfysize=2in \epsffile{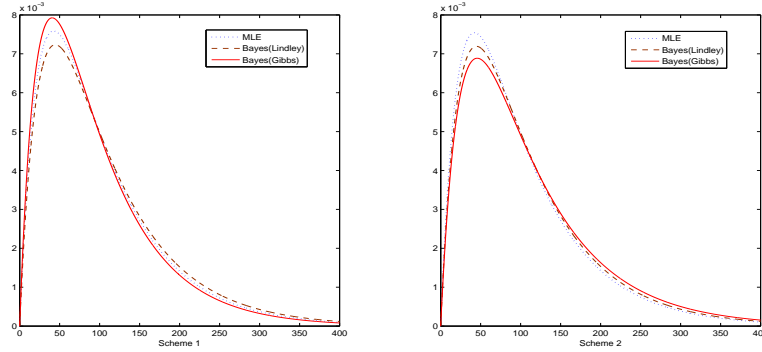}}
\vspace{-0.3in}
\caption{\scriptsize Estimated density functions with informative priors for scheme 1(left) and scheme 2(right). }
\end{figure}
 \label{sch}\input{epsf}
\epsfxsize=3in \epsfysize=2in
\begin{figure}
\centerline{\epsfxsize=5in \epsfysize=2in \epsffile{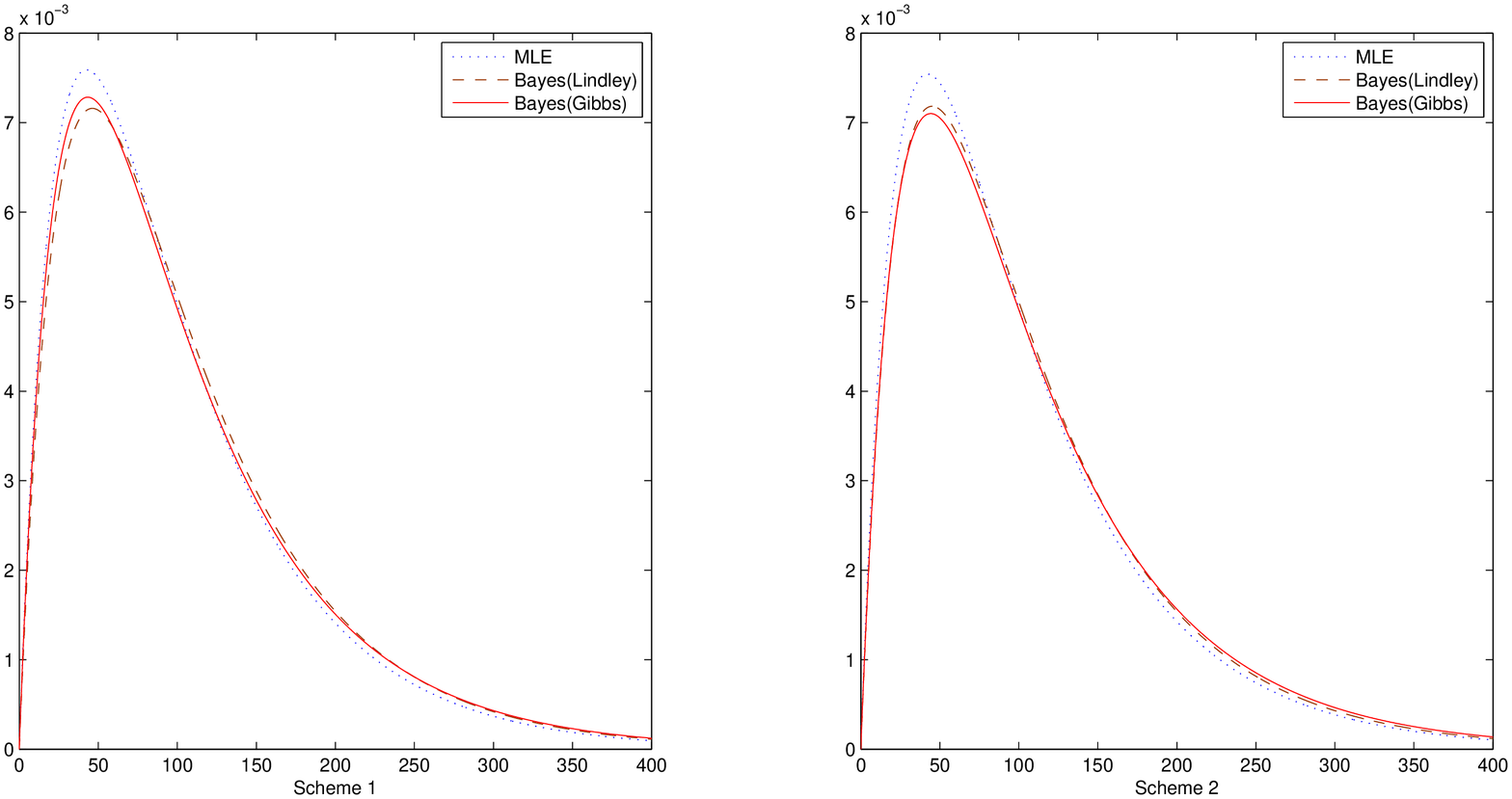}}
\vspace{-0.3in}
\caption{\scriptsize Estimated density functions with non-informative priors for scheme 1(left) and scheme 2(right). }
\end{figure}
Comparing the two schemes with informative and non-informative priors, it is observed that for scheme 1, estimators have smaller standard errors than scheme 2, as expected. Also it is clear that the Bayes estimators depend on the hyper parameters. Because the HPD confidence intervals based on informative priors are slightly smaller than corresponding length of HPD confidence intervals based on non-informative priors, therefore the prior informative should be used if they are available.

\section{Conclusion}
\par In this article, we have studied the classical and Bayes inference procedure for the Type-II hybrid censored $\mbox{WE}(\alpha,\lambda)$ distribution. We provide the maximum likelihood estimators and it is observed that the maximum likelihood estimators of the unknown parameters can not be obtained in the closed form and we suggest the EM algorithm to compute them. we also earn the Bayes estimators of the unknown parameters and show that they can not be obtained in explicit forms, and we have proposed two approximation methods to earn them. We have compared the performance of the different methods by Monte Carlo simulations, and it is observed that the performance of quite satisfactory.

\newpage
\begin{center}
{\bf Appendix A}
\end{center}
{{\label{ap}\theorem{Given $X_{(1)}=x_{(1)},\cdots,X_{(r)}=x_{(r)}$, the conditional distribution of $Z_i$ for $i=1,\cdots,n-r$ is
$$f_{Z|X}(z_i|X_{(1)}=x_{(1)},\cdots,X_{(i)}=x_{(i)})=f_{Z|X}(z_i|X_{(i)}=x_{(i)})=
\frac{f_{WE}(z_i)}{1-F_{WE}(x_i)}$$
where\\ $$f_{WE}(x)=\frac{\alpha+1}{\alpha^2}\beta e^{-\frac{\beta}{\alpha}x}(1-e^{-\beta x})$$
and $$F_{WE}(x)=1-\frac{1}{\alpha}e^{-\frac{\beta}{\alpha}x}[\alpha+1-e^{-\beta x}]$$}}}
{\proof The proof can be obtained similarly as in Ng et al. (2002).}\\
Note that using Theorem \ref{ap}, we can write
\begin{equation*}
\displaystyle
\begin{array}{c}
A(c;\alpha,\beta)=E[Z_i|Z_i>c]=\int_c^\infty \frac{\alpha+1}{\alpha}\frac{\beta}{K}xe^{-\frac{\beta}{\alpha}x}(1-e^{-\beta x})dx\;\;\;\;\left(K=e^{-\frac{\beta}{\alpha}c}(\alpha+1-e^{-\beta c})\right)\\ \\
\hspace{1in}=\frac{(\alpha+1)\alpha}{K\beta}\int_{\frac{\beta}{\alpha}c}^\infty ue^{-u}du-\frac{\alpha}{K\beta(\alpha+1)}\int_{\frac{\beta(\alpha+1)}{\alpha}c}^\infty ve^{-v}dv \;\;\;\left(put \;\;u=\frac{\beta}{\alpha}x\;\;\&\;\;v=\frac{\beta(\alpha+1)}{\alpha}x\right)\\ \\
\hspace{-0.5in}=\frac{(\alpha+1)\alpha}{K\beta}(-e^{-u})(u+1)]_{\frac{\beta}{\alpha}c}^\infty-\frac{\alpha}{K\beta(\alpha+1)}(-e^{-v})(v+1)]
_{\frac{\beta(\alpha+1)}{\alpha}c}^\infty\\ \\
\hspace{-0.7in}=\frac{(\alpha+1)\alpha}{K\beta}(\frac{\beta}{\alpha}c+1)e^{-\frac{\beta}{\alpha}c}-\frac{\alpha}{K\beta(\alpha+1)}(\frac{\beta(\alpha+1)}{\alpha}c+1)e^{
-\frac{\beta(\alpha+1)}{\alpha}c}\\ \\
\hspace{-1.25in}=\frac{\alpha e^{-\frac{\beta}{\alpha}c}}{K\beta}\left[(\alpha+1)(\frac{\beta}{\alpha}c+1)-\frac{e^{-\beta c}}{\alpha+1}(\frac{\beta(\alpha+1)c}{\alpha}+1)\right]\\ \\
\hspace{-1.1in}=\frac{e^{-\frac{\beta}{\alpha}c}}{K\beta}\left[(\alpha+1)(\beta c+\alpha)-\frac{e^{-\beta c}}{\alpha+1}\left(\beta(\alpha+1)c+\alpha\right)\right]
\\ \\ \displaystyle
\hspace{-1.9in}=\frac{(\alpha+1)(\beta c+\alpha)-\frac{e^{-\beta c}(\beta c(\alpha+1)+\alpha)}{\alpha+1}}{\beta(\alpha+1-e^{-\frac{\beta}{\alpha}c})}
\end{array}
\end{equation*}
and about $B(c;\alpha,\beta)$, we have:
$$B(c;\alpha,\beta)=E[\ln (1-e^{-\beta Z_i})|Z_i>c]=\int_c^\infty\frac{\alpha+1}{\alpha}\frac{\beta}{K}\ln(1-e^{-\beta x})e^{-\frac{\beta}{\alpha}x}(1-e^{-\beta x})dx$$$$
=-\frac{\alpha+1}{\alpha K}\int_{1-e^{-\beta c}}^1u(1-u)^{\frac{1}{\alpha}-1}\ln(u)$$
$$=\frac{\alpha+1}{\alpha K}\left[\frac{1}{18}(1-e^{-\beta c})^3(-\frac{1}{\alpha}+1)\mbox{hypergeom}_{3,2}([3,3,-\frac{1}{\alpha}+2],[4,4],1-e^{-\beta c})\right.$$$$+\left(-\frac{1}{4}\mbox{hypergeom}_{2,1}([2,-\frac{1}{\alpha}+1],[3],1-e^{-\beta c})\right.$$$$\left.+\frac{1}{2}\ln(x)\mbox{hypergeom}_{2,1}([2,-\frac{1}{\alpha}+1],[3],1-e^{-\beta c})\right)(1-e^{-\beta c})^2$$$$-\frac{1}{18}(-\frac{1}{\alpha}+1)\mbox{hypergeom}_{3,2}([3,3,-\frac{1}{\alpha}+2],[4,4],1)$$$$\left.+\frac{1}{4}
\mbox{hypergeom}_{2,1}([2,-\frac{1}{\alpha}+1],[3],1)\right]$$
where hypergeom(.) is Generalized hypergeometric function. This function is also known as Barnes's extended hypergeometric function. The definition of $F_{p,q}(\mathbf{n},\mathbf{d},\xi)$ is:
$$F_{p,q}(\mathbf{n},\mathbf{d},\xi)=\sum_{k=0}^\infty\frac{\lambda^k\Pi_{i=1}^p\Gamma(n_i+k)\Gamma^{-1}(n_i)}{\Gamma(k+1)
\Pi_{i=1}^q\Gamma(d_i+k)\Gamma^{-1}(d_i)},$$
where $\mathbf{n}=[n_1,\cdots,n_p]$, $p$ is the number of operands of $\mathbf{n}$, $\mathbf{d}=[d_1,\cdots,d_q]$ and $q$ is the number of operands
of $\mathbf{d}$. Generalized hypergeometric function is quickly evaluated and readily available in standard software such as Maple.

\begin{center}
{\bf Appendix B}
\end{center}
The conditional density of $\beta$ given $\alpha$ and $\b{x}$ is
$$\pi(\beta|\alpha,\b{x})\propto\beta^{w_2+r-1}e^{-\beta \left(w_1+\frac{\beta}{\alpha}\sum_{i=1}^rx_i+c(n-r)\right)}
(\alpha+1-e^{-\beta c})^{n-r}\Pi_{i=1}^r(1-e^{-\beta x_i}).$$
This function is log-concave because we have
$$\frac{\partial^2\log(\pi(\beta|\alpha,\b{x}))}{\partial\beta^2}=-\frac{w_2+r-1}{\beta^2}-\frac{(n-r)c^2e^{-\beta c}(\alpha+1)}{(\alpha+1-e^{-\beta c})^2}-\sum_{i=1}^r\frac{x_i^2e^{-\beta x_i}}{(1-e^{-\beta x_i})^2}<0.$$
Therefore, the result follows.
\\ \quad\\

\begin{center}
{\bf Appendix C}
\end{center}
The conditional distribution of $\alpha$ given $\beta$ and $\b{x}$ is
$$\pi(\alpha|\beta,\b{x})\propto\alpha^{w_4-n-r-1}(\alpha+1)^re^{-\alpha w_3}e^{-\frac{\beta}{\alpha}(\sum_{i=1}^rx_i+c(n-r))}
(\alpha+1-e^{-\beta c})^{n-r}.$$
In this function, we have $\pi(\infty|\beta,\b{x})\rightarrow\pi(0|\beta,\b{x})=0$ and $\pi(\alpha|\beta,\b{x})\geq0~~\forall\alpha$, now it is enough that prove $\pi(\alpha|\beta,\b{x})$ is bounded. With simple calculation we see that $\forall\alpha$ this function is less than the gamma function and the gamma function is a bounded function, so this function is bounded. Therefore $\pi(\alpha|\beta,\b{x})$ has a finite maximum point.

\bibliographystyle{plain}

\end{document}